\documentclass[a4paper,11pt]{article}
\usepackage[T1]{fontenc}
\usepackage{lmodern,amsmath,amssymb,amsthm}
\usepackage[margin=26mm,top=22mm,bottom=22mm]{geometry}
\usepackage{microtype}
\usepackage{xcolor}
\usepackage[bottom]{footmisc}
\usepackage[bookmarksnumbered=true]{hyperref}
\definecolor{linkink}{RGB}{29,66,89}
\hypersetup{colorlinks=true,linkcolor=linkink,citecolor=linkink,urlcolor=linkink,
 pdftitle={Frechet-Urysohn Scott spaces and criteria for sobriety},
 pdfsubject={Domain theory and non-Hausdorff topology}}
\newcommand{\N}{\mathbb N}
\newcommand{\up}{\mathord{\uparrow}}
\newcommand{\dn}{\mathord{\downarrow}}
\newcommand{\cl}{\operatorname{cl}}
\newcommand{\Int}{\operatorname{int}}

\newtheorem{theorem}{Theorem}[section]
\newtheorem{lemma}[theorem]{Lemma}
\newtheorem{proposition}[theorem]{Proposition}
\newtheorem{corollary}[theorem]{Corollary}
\theoremstyle{definition}
\newtheorem{definition}[theorem]{Definition}

\theoremstyle{remark}

\numberwithin{equation}{section}
\title{\bfseries The Scott sobriety on countable Fr\'echet-Urysohn dcpos and locally compact meet-continuous dcpos %
  \thanks{This work is supported by the National Natural Science Foundation of China (Grant nos.12331016,12471438) and the Sichuan Science and Technology Program (Grant No. 2026NSFSC0784).}}

\author{Zhengmao He\thanks{Email: \href{hezhengmaomath@163.com}{\nolinkurl{hezhengmaomath@163.com}}.}\\[0.4em]
  {\small School of Sciences, Southwest Petroleum University, Chengdu, 610500, China}}

\date{}
\begin{document}
\newlength{\abstractlabelwidth}
\renewenvironment{abstract}
  {%
    \par\noindent
    \settowidth{\abstractlabelwidth}{\textit{Keywords.}}%
    \makebox[\abstractlabelwidth][r]{\textit{Abstract.}}%
    \space\small\ignorespaces
  }
  {\par\smallskip}
\maketitle
%\address{School of Sciences, Southwest Petroleum University,Chengdu 610500,Sichuan China}
\begin{abstract}
In this paper, we focus on the Scott Sobriety of countable dcpos and locally compact dcpos. Recently, Chen, Kou and Yang proved that every countable dcpo with its Scott topology first-countable is sober endowed with the Scott topology. A natural question is whether a countable dcpo whose Scott topology is a Fr\'echet-Urysohn sapce must be sober. We show that a countable dcpo with condition FG is Scott sober if its Scott topology is Fr\'echet-Urysohn. Besides, a countable Scott non sober Fr\'echet-Urysohn dcpo is constructed. He and Zhang gave two sufficient conditions for the open problem whether every locally compact dcpo is Scott sober. One of their sufficient conditions is that any locally compact meet-continuous coherent dcpo is Scott sober. We further prove that each locally compact meet-continuous dcpo is always sober with respect to the Scott topology.
\end{abstract}


%\blfootnote{This work is supported by the National Natural Science Foundation of China (Grant nos.12471438,12601906) and the Sichuan Science and Technology Program (Grant No. 2026NSFSC0784).}
%\blfootnote{E-mail address: hezhengmaomath@163.com (Z.M.He).}

\section{Introduction}\label{sec:intro}

Sobriety is probably the most important and useful property of non-Hausdorff topological spaces. In Domain theory, a classical result is that every quasicontinuous dcpos and in particular
every continuous dcpos is sober with respect to the Scott topology (\cite{domains}). Countability alone gives some different examples. Johnstone constructed a countable dcpo with
a non-sober Scott space (\cite{johnstone}). Miao, Xi, Li and Zhao later
answered a question of Jung by constructing a countable complete
lattice with non-sober Scott topology. Their construction even admits
a distributive lattice order (\cite{miao2023}).

Countability nevertheless has a strong consequence for convergence.
The Scott space of every countable dcpo is sequential. Indeed, a
sequentially closed subset is lower, as the constant sequence at $p$
converges to each point below $p$. It also contains the supremum of
each of its directed subsets, since every countable directed set has
an increasing cofinal sequence converging to its supremum. Chen, Kou
and Yang obtained a positive sobriety result under a stronger assumption.
They proved that first countability of $\Sigma P$ forces a countable
dcpo $P$ to be quasicontinuous and hence Scott sober
(\cite{Chenkouyang}). These two results single out the
Fr\'echet--Urysohn property. First countability implies this property,
which in turn implies sequentiality. Its distinctive requirement is
that each individual point of $\cl A$ be the limit of a sequence
from $A$. A related theorem of He concerns the square of a space.
An $\omega$-well-filtered coherent $d$-space $X$ is sober if
$X\times X$ is Fr\'echet--Urysohn (\cite{hefu}). Here we ask whether the
Fr\'echet-Urysohn property of $\Sigma P$ itself suffices when $P$
is countable.

We answer this question negatively. The countable dcpo constructed in
Section 4 has an irreducible Fr\'echet-Urysohn Scott
space and no greatest element. Its whole space is therefore a closed
irreducible set without a generic point. To verify the
Fr\'echet-Urysohn property we compute the Scott closure of an arbitrary
subset. One sequence in the construction, $(a^1_r)_{r\ge1}$,
converges to every point of the space.

The positive result depends on the order structure of the common upper
bounds. Suppose that for every $x,y\in P$, there is a finite set $F_{x,y}\subseteq P$ such that
\[
                     \up x\cap\up y=\up F_{x,y}.
\]
Miao, Yuan and Li used this finite generation property in their study
of join continuity and well-filteredness (\cite{wellfiltered}). Here
its role is to recover directed approximation from a single convergent
sequence. Given $a_n\to x$, countability allows us to pass to a
subsequence on which, for each $p\in P$, the assertion $p\le a_n$
has a constant truth value from some index onwards. Finite generation
makes the set of eventual lower bounds directed. This yields weak one-step closure in
the sense of Miao, Li and Zhao (\cite{onestep}). We then prove that
every countable dcpo with weak one-step closure is Scott sober. In that proof, an enumeration of an irreducible
closed set without a greatest element would yield an increasing branch
with no upper bound in the set, contrary to Scott closedness.
Thus the Fr\'echet-Urysohn property suffices
under FG. In particular, the counterexample constructed by Miao,Xi, Li and Zhao
cannot be Fr\'echet-Urysohn, since binary joins give
$\up x\cap\up y=\up(x\vee y)$(\cite{miao2023}).

A second question arises when the cardinality restriction is lifted.
Xu and Zhao asked whether local compactness of $\Sigma P$ implies
sobriety for an arbitrary dcpo $P$ (\cite{xuzhao}).
The countable case has an affirmative answer. In fact, He, Dong and
Wang proved that a countable dcpo with core-compact Scott space is
quasicontinuous (\cite{hedongwang}). For arbitrary
cardinality, He and Zhang proved sobriety when $\Sigma P$ is locally
compact and second countable. They also proved it when $P$ is
meet-continuous and $\Sigma P$ is locally compact and coherent
(\cite{hz}).

{\bf Is every meet-continuous locally compact dcpo Scott sober? }\\
Chen constructed a meet-continuous
dcpo with a compact second-countable core-compact dcpo that
is neither locally compact nor sober (\cite{chencore}).
This gives a negative answer to the related question of Jia recorded in
(\cite{xuzhao}). Theorem 5.18 answers Question~3.14 affirmatively for
meet-continuous dcpos of arbitrary cardinality. He and Zhang use coherence
to make $\up a\cap\up b$ compact(\cite{hz}). This allows them to take infima
along the decreasing columns of a triangular array. We instead work
inside each of a sequence of prescribed Scott-compact sets $K_n$.
Meet-continuity supplies common lower bounds for finite families in
suitably chosen neighbourhoods. By successive relatively closed
refinements and a finitely branching tree argument, we obtain points
$p_n\in K_n$ with $p_n\le p_{n+1}$ whenever
$K_{n+1}\subseteq\Int(\up K_n)$ for all $n$. Inside an irreducible
Scott-closed set $F$, this selection puts a point in every countable
intersection of nonempty relatively open sets. We can then adapt the
minimal closed set argument of He and Zhang to show that the points
with a countable neighbourhood base in $F$ form a dense directed set.
Its supremum is the greatest element of $F$.

\section{Preliminaries}\label{c:prelim}

Write\(\N_+=\{1,2,\ldots\}\) and \(\N_0=\{0,1,2,\ldots\}\). A countable
set is a finite or countably infinite set.
For a subset \(A\) of a partially ordered set \(P\), put
\[
  \up A=\{y\in P\mid\exists a\in A, a\le y\},
  \qquad
  \dn A=\{y\in P\mid\exists a\in A, y\le a\}.
\]
For a singleton, we write \(\up x\) and \(\dn x\). A set \(A\) is \emph{upper} if
\(A=\up A\) and \emph{lower} if \(A=\dn A\).

A subset \(D\) of \(P\) is \emph{directed} if it is nonempty and every
pair of its elements has an upper bound in \(D\). A \emph{dcpo} is a
partially ordered set in which every directed subset \(D\) admits a
supremum, denoted by \(\bigvee D\). A subset \(E\subseteq D\) is
\emph{cofinal} in \(D\) if each \(d\in D\) lies below some \(e\in E\).A dcpo $P$ is {\em meet-continuous} if for every directed subset $D$ and every $x\in P$, $x\leq\bigvee D$ implies $x\in\overline{(\dn x\cap\dn D)}_{\Sigma P}$.

Let $P$ be a poset. For each pair of $x, y\in P$, we say that $x$ is {\em way below} $y$, in symbols $x\ll y$, if for each directed set $D\subseteq P$ whenever $\bigvee D$ exists, $y\leq\bigvee D$ implies $\up x\cap D\neq\emptyset$. A poset $P$ is {\em continuous} if $\Downarrow p=\{b\in P\mid b\ll p\}$ is directed and $\bigvee \Downarrow p=p$ for every $p\in P$.

Let \(P\) be a dcpo. Its \emph{Scott topology} \(\sigma(P)\) consists of
the upper sets \(U\) with the following property. For each directed
\(D\subseteq P\), if \(\bigvee D\in U\), then \(D\cap U\ne\varnothing\).
We denote the space \((P,\sigma(P))\) by \(\Sigma P\). The Scott closed sets are precisely the lower sets containing the suprema
of their directed subsets.

A nonempty subset \(A\) of a space \(X\) is \emph{irreducible} if \(A\subseteq C_1\cup C_2\) with \(C_1,C_2\) closed in \(X\)
implies \(A\subseteq C_1\) or \(A\subseteq C_2\). Equivalently, if
two open sets meet \(A\), then their intersection meets \(A\). A \(T_0\) space
is \emph{sober} if each nonempty irreducible closed set is the closure of
a unique point.

A subset $A$ of a topological space $X$ is {\em saturated} if $A=\up A$. Given a topological space $X$, we shall use $K(X)$ to denote the poset of all nonempty
compact saturated subsets of $X$ with the reverse inclusion order.

A topological space $X$ is called {\em locally compact} if for each $U\in \mathcal{O}(X)$ and for each $x\in U$, there exist an open set $V$ and a $K\in K(X)$ such that $x\in V\subseteq K\subseteq U$. A topological space $X$ is said to be  {\em core-compact } if  $\mathcal{O}(X)$ is continuous.

A space \(X\) is \emph{Fr\'echet--Urysohn} if for each \(A\subseteq X\)
and \(x\in\cl_X A\) there exists a sequence of points of \(A\) converging
to \(x\). If \(x_n\to x\) in
\(\Sigma P\), then \(x_n\to y\) for every \(y\le x\), since every open
neighbourhood of \(y\) contains \(x\). If \(x_n\to x\) and \(x_n\le y_n\) for all \(n\),
then also \(y_n\to x\), since any open set containing \(x_n\) contains
\(y_n\).

In a topological space $X$, a \(G_\delta\) set is a countable intersection
of open sets and a subset is dense precisely when it meets every
nonempty open set.

\begin{definition}[\cite{onestep}]\label{c:wosc}
A dcpo \(P\) has \emph{weak one-step closure} if for every \(A\subseteq P\)
and \(x\in\cl A\), there is a directed \(D\subseteq\dn A\) with
\(x\le\bigvee D\).
\end{definition}

For any dcpo, if \(D\subseteq\dn A\) is directed, then \(\cl A\) contains
\(D\), its supremum and all points below that supremum. Thus the definition
is equivalent to the identity
\[
     \cl A
       =\{x\in P:x\le\bigvee D\text{ for some directed }D\subseteq\dn A\}
       (A\subseteq P).
\]

\begin{definition}\label{c:fg}
A poset \(P\) satisfies the \emph{finite generation condition for
intersections of principal upper sets} if for each $x,y\in P$,
\begin{equation}\tag{FG}\label{c:eqfg}
  \up x\cap\up y=\up F_{x,y}\text{ for some finite }F_{x,y}\subseteq P.
\end{equation}
\end{definition}

\begin{lemma}\label{c:closed}
Let $F$ be a Scott closed subset of a dcpo $P$. If $P$ has weak
one-step closure, then $F$ has weak one-step closure.
\end{lemma}

\begin{proof}
For every directed set $D\subseteq F$, Scott closedness gives
$\bigvee_P D\in F$. Since $F$ carries the induced order,
\[
\bigvee_F D=\bigvee_P D.
\]
Thus $F$ is a dcpo. Let $A\subseteq F$ and $x\in\operatorname{cl}_{\Sigma F}A$.
Put $C=\operatorname{cl}_{\Sigma P}A$. Since $A\subseteq F$ and
$F$ is Scott closed in $P$, we have $C\subseteq F$. Moreover, $C$
is lower in $F$, and every directed $D\subseteq C$ satisfies
$\bigvee_F D=\bigvee_P D\in C$. Hence $C$ is Scott closed in $F$.
So
\[
x\in\operatorname{cl}_{\Sigma F}A
\subseteq C=\operatorname{cl}_{\Sigma P}A.
\]

By weak one-step closure of $P$, there exists a directed set
$D\subseteq\dn_{P} A$ such that $x\leq\bigvee_P D$.
Since $A\subseteq F$ and $F$ is lower in $P$,
\[
\dn_P A\subseteq F\mbox{ and }
\dn_F A=F\cap\dn_P A=\dn_P A.
\]
Consequently, $D\subseteq\dn_F A$ and
$x\leq\bigvee_P D=\bigvee_F D$. This proves that $F$ has weak one-step closure.
\end{proof}

\section{A sufficient condition for Fr\'echet-Urysohn countable dcpo becoming sober}\label{c:positive}

In this section, we will give a sufficient condition that ensures a countable Fr\'echet-Urysohn dcpo becoming Scott sober

\begin{lemma}\label{lem:infinitely-many-lower-bound}
{\rm Let $P$ be a dcpo and $x_{n}$ a sequence converging to $x$ in $\Sigma P$. Then for every
Scott open neighbourhood $U$ of $x$, there exists $u\in U$ such that
$u\leq x_n$ for infinitely many $n$.}
\end{lemma}

\begin{proof}
Fix $x\in U\in\sigma(P)$. Suppose that for every $u\in U$, the set
$I(u)=\{n\in\mathbb N_+:u\leq x_n\}$ is finite. Since $x\in U$, we can choose $N\geq1$ such that
$x\not\leq x_n$ for every $n\geq N$. Put
\[
C=(P\setminus U)\cup\bigcup_{n\geq N}\downarrow x_n.
\]
Then $C$ is lower, $x\notin C$ and $x_n\in C$ for every $n\geq N$. Let $D\subseteq C$ be a directed subset and put $d=\bigvee D$. If $d\notin U$,
then $d\in C$. Suppose $d\in U$. As $U$ is Scott open, there is a
$d_0\in D\cap U$ and put $D_0=D\cap\uparrow d_0$. By directedness of
$D$, $D_0$ is directed and cofinal in $D$. Hence $\bigvee D_0=d$.
For every $z\in D_0$, we have $z\in C\cap U$. So there is a
$N\leq n$ such that $d_0\leq z\leq x_n$. Consequently,
$$D_0\subseteq
\bigcup_{\substack{n\in I(d_0)\\ n\geq N}}\dn x_n
\subseteq C.
$$
Since $d_0\in U$, the middle union above is a finite union. It therefore contains $\bigvee D_0=d$. Hence $C$ is Scott closed. Now $P\setminus C$ is a Scott open neighbourhood of
$x$. But $x_n\notin P\setminus C$ for every $n\geq N$, a contradiction.
\end{proof}

\begin{lemma}\label{lem:eventual-lower-bounds}
{\rm Let $P$ be a dcpo and $x_{n}$ a sequence converging to $x$ in $\Sigma P$. Then there is
a subsequence $(x_{n_k})_{k\geq1}$ such that for each $p\in P$ the
value of $p\leq x_{n_k}$ is constant for all sufficiently large
$k$. Its set of eventual lower bounds
\[
L=\{p\in P:\exists k_0\geq1,\forall k\geq k_0,\ p\leq x_{n_k}\}
\]
is lower and satisfies $x\in\operatorname{cl}_{\Sigma P}L$.}
\end{lemma}

\begin{proof}
Write $P=\{p_1,p_2,\ldots\}$, repeating elements if $P$ is finite,
and put $I_0=\mathbb N_+$. Given an infinite $I_{i-1}$, at least one
of the sets
\[
\{n\in I_{i-1}:p_i\leq x_n\},
\qquad
\{n\in I_{i-1}:p_i\not\leq x_n\}
\]
is infinite. Choose such a set as $I_i$. Choose successively
$n_1<n_2<\cdots$ with $n_k\in I_k$. For every $k\geq i$,
\[
n_k\in I_k\subseteq I_i.
\]
Thus either $p_i\leq x_{n_k}$ for every $k\geq i$, or
$p_i\not\leq x_{n_k}$ for every $k\geq i$.

If $q\leq p\in L$, then $q\leq p\leq x_{n_k}$ for all sufficiently
large $k$, so $q\in L$. Hence $L$ is lower. For any $x\in U\in\sigma(P)$,
the subsequence still converges to $x$. By
Lemma~\ref{lem:infinitely-many-lower-bound}, there exists $p\in U$
with $p\leq x_{n_k}$ for infinitely many $k$. By the eventual constancy
just proved,
\[
(\exists k_0\geq1)(\forall k\geq k_0)\quad p\leq x_{n_k}.
\]
Thus $p\in U\cap L$. Since every Scott-open neighbourhood of $x$
meets $L$, we have $x\in\operatorname{cl}_{\Sigma P}L$.
\end{proof}

\begin{proposition}
Let $P$ be a countable dcpo satisfying \textup{(FG)}. If $\Sigma P$
is Fr\'echet--Urysohn, then $P$ has weak one-step closure.
\end{proposition}

\begin{proof}
Let $A\subseteq P$ and $x\in\operatorname{cl}_{\Sigma P}A$.
By the Fr\'echet--Urysohn property, there is a sequence in $A$
converging to $x$. Pass to the subsequence given by
Lemma~\ref{lem:eventual-lower-bounds}, denote it by $(a_n)_{n\geq1}$,
and let $L$ be its set of eventual lower bounds. Then
\[
x\in\operatorname{cl}_{\Sigma P}L,\qquad
\varnothing\neq L\subseteq\downarrow A.
\]

Fix $p,q\in L$. Choose $N\geq1$ such that $p,q\leq a_n$ for all
$n\geq N$. By \textup{(FG)}, there is a finite $F_{p,q}\subseteq P$
such that
\[
a_n\in\uparrow p\cap\uparrow q
     =\uparrow F_{p,q}\qquad(n\geq N).
\]
Hence
\[
\{n\in\mathbb N_+:n\geq N\}
=
\bigcup_{f\in F_{p,q}}\{n\geq N:f\leq a_n\}.
\]
Since $F_{p,q}$ is finite, some $f\in F_{p,q}$ satisfies $f\leq a_n$
for infinitely many $n$. Eventual constancy gives $f\leq a_n$ for
all sufficiently large $n$, so $f\in L$. Moreover,
\[
f\in F_{p,q}\subseteq\uparrow F_{p,q}
 =\uparrow p\cap\uparrow q.
\]
Thus $p,q\leq f$, proving that $L$ is directed.

Put $s=\bigvee L$. Since $\downarrow s$ is Scott closed and contains
$L$, we have $\operatorname{cl}_{\Sigma P}L\subseteq\downarrow s$.
Conversely, $\operatorname{cl}_{\Sigma P}L$ contains the directed set
$L$, hence contains $s$ and every point below $s$. Therefore
\[
x\in\operatorname{cl}_{\Sigma P}L=\downarrow s,
\qquad x\leq\bigvee L.
\]
Taking $D=L\subseteq\downarrow A$ proves weak one-step closure.
\end{proof}

\begin{lemma}\label{lem:weak-one-step-FU}
If a countable dcpo $P$ has weak one-step closure, then $\Sigma P$
is Fr\'echet--Urysohn.
\end{lemma}

\begin{proof}
Let $A\subseteq P$ and $x\in\operatorname{cl}_{\Sigma P}A$.
Choose a directed set $D\subseteq\downarrow A$ with $x\leq\bigvee D$.
Write $D=\{e_1,e_2,\ldots\}$, repeating elements if $D$ is finite.
Put $d_1=e_1$. For each $n\geq1$, choose $d_{n+1}\in D$ with
$d_n\leq d_{n+1}$ and $e_{n+1}\leq d_{n+1}$. Then $(d_n)_{n\geq1}$
is increasing and cofinal in $D$, so $\bigvee_{n\geq1}d_n=\bigvee D$.

Since $d_n\in\downarrow A$, choose $a_n\in A$ with $d_n\leq a_n$.
Let $U\in\sigma(P)$ contain $x$. By upperness, $\bigvee D\in U$.
Scott openness gives $d_N\in U$ for some $N$. For every $n\geq N$,
we have $d_N\leq d_n\leq a_n$, whence $a_n\in U$. Thus $a_n\to x$.
\end{proof}

\begin{lemma}\label{lem:finite-upper-interior}
Let $P$ be a countable dcpo with weak one-step closure. Suppose that
$\Sigma P$ is irreducible and $P$ has no greatest element. For every
nonempty $U\in\sigma(P)$, there exists a finite set $E\subseteq U$
such that $\operatorname{int}_{\Sigma P}(\uparrow E)\neq\varnothing$.
\end{lemma}

\begin{proof}
Let $M$ be the set of maximal elements of $P$. Fix $x\in P$.
Every nonempty chain in $\uparrow x$ has its supremum in $\uparrow x$,
and $x$ is an upper bound of the empty chain. By Zorn's lemma,
$x\leq m$ for some $m\in M$. Hence $P=\downarrow M$, and therefore
$\operatorname{cl}_{\Sigma P}M=P$.

For each $m\in M$, the equality $M=(M\setminus\{m\})\cup\{m\}$ gives
$$P=\operatorname{cl}_{\Sigma P}(M\setminus\{m\})\cup\downarrow m.$$
Both sets on the right are Scott closed. Since $P$ has no greatest
element, $\downarrow m\neq P$. Irreducibility now gives
$\operatorname{cl}_{\Sigma P}(M\setminus\{m\})=P$.

Fix a nonempty $U\in\sigma(P)$ and choose $m\in M\cap U$.
For every Scott-open neighbourhood $V$ of $m$, the nonempty open set
$V\cap U$ meets $M\setminus\{m\}$. Thus
$m\in\operatorname{cl}_{\Sigma P}((M\cap U)\setminus\{m\})$.
By Lemma~\ref{lem:weak-one-step-FU}, there is a sequence in
$(M\cap U)\setminus\{m\}$ converging to $m$. Each value $v$ occurs
only finitely often, since $m\not\leq v$ and $P\setminus\downarrow v$
is an open neighbourhood of $m$. Passing to a subsequence, choose
pairwise distinct $m_n\in(M\cap U)\setminus\{m\}$ with $m_n\to m$.

Suppose that $\operatorname{int}_{\Sigma P}(\uparrow E)=\varnothing$
for every finite $E\subseteq U$. Enumerate $U=\{u_1,u_2,\ldots\}$.
Put $E_n=\{u_1,\ldots,u_n\}$ and $A_n=P\setminus\uparrow E_n$.
Then $A_n$ is lower and
$$\operatorname{cl}_{\Sigma P}A_n
=P\setminus\operatorname{int}_{\Sigma P}(\uparrow E_n)=P.$$
By weak one-step closure, there exists a directed set
$D'_n\subseteq\downarrow A_n=A_n$ with $m_n\leq\bigvee D'_n$.
Maximality gives $\bigvee D'_n=m_n$. Since $m_n\in U$ and $U$ is
Scott open, choose $d_0\in D'_n\cap U$. Put $D_n=D'_n\cap U$.
Each $d\in D'_n$ has a common upper bound with $d_0$ in $D'_n$,
and every such bound belongs to $U$. Thus $D_n$ is cofinal in $D'_n$.
It is also directed, since $U$ is upper. Consequently,
$D_n\subseteq U\setminus\uparrow E_n$ and $\bigvee D_n=m_n$.

Put $B=\bigcup_{n\geq1}D_n$. The Scott-closed set
$\operatorname{cl}_{\Sigma P}B$ contains $\bigvee D_n=m_n$ for every
$n$, and hence contains their limit $m$. By
Lemma~\ref{lem:weak-one-step-FU}, choose a sequence $(b_k)_{k\geq1}$
in $B$ converging to $m$. For each fixed $n$, we have
$D_n\subseteq\downarrow m_n$ and $m\not\leq m_n$. Since
$P\setminus\downarrow m_n$ is an open neighbourhood of $m$, the set
$\{k:b_k\in D_n\}$ is finite.

Fix $i\geq1$. For $n\geq i$, we have $u_i\in E_n$, so
$D_n\cap\uparrow u_i=\varnothing$. As each $b_k$ belongs to $B$,
it follows that
$$\{k:u_i\leq b_k\}\subseteq\bigcup_{1\leq n<i}\{k:b_k\in D_n\}.$$
The right-hand side is finite. Thus no element of $U$ lies below
infinitely many $b_k$, contrary to Lemma~3.1 applied to $b_k\to m$
and the neighbourhood $U$ of $m$.
\end{proof}

\begin{theorem}\label{thm:countable-weak-one-step-sober}
Every countable dcpo with weak one-step closure has a sober Scott space.
\end{theorem}

\begin{proof}
Let $P$ be a countable dcpo with weak one-step closure, and let $F$
be a nonempty irreducible Scott-closed subset of $P$. By Lemma~2.3,
$F$ is a countable dcpo with weak one-step closure. Directed suprema
in $F$ agree with those in $P$. Since $F$ is lower in $P$, every
Scott-closed subset of $F$ is Scott closed in $P$. Conversely, the
intersection of $F$ with any Scott-closed subset of $P$ is Scott
closed in $F$. Thus $\Sigma F$ carries the subspace topology and is
irreducible.

Suppose that $F$ has no greatest element. Enumerate
$F=\{p_1,p_2,\ldots\}$, allowing repetitions, and put $W_0=F$.
Given a nonempty $W_{n-1}\in\sigma(F)$, the set
$F\setminus\downarrow_F p_n$ is nonempty and Scott open in $F$.
By irreducibility, $U_n=W_{n-1}\setminus\downarrow_F p_n$ is
nonempty. Apply Lemma~\ref{lem:finite-upper-interior} in $F$ to
choose a finite $E_n\subseteq U_n$ such that
$W_n=\operatorname{int}_{\Sigma F}(\uparrow_F E_n)\neq\varnothing$.
Since $U_n$ is upper, $W_n\subseteq\uparrow_F E_n\subseteq U_n$.
In particular, $E_n\neq\varnothing$ and
$E_n\cap\downarrow_F p_n=\varnothing$. Moreover,
$$E_{n+1}\subseteq U_{n+1}\subseteq W_n\subseteq\uparrow_F E_n.$$

Consider the tree of finite increasing lists $(e_1,\ldots,e_n)$ with
$e_i\in E_i$, ordered by extension and rooted at the empty list.
Each level is nonempty. Indeed, choose any $e_n\in E_n$ and use
$E_{i+1}\subseteq\uparrow_F E_i$ successively to choose
$e_i\in E_i$ with $e_i\leq e_{i+1}$ for $i=n-1,\ldots,1$.
Each node at level $n$ has at most $|E_{n+1}|$ children. By
K\"onig's lemma, the tree has an infinite branch $(e_n)_{n\geq1}$
with $e_n\in E_n$ and $e_n\leq e_{n+1}$ for every $n$.
Its supremum $p=\bigvee_{n\geq1}e_n$ belongs to $F$. Choose $N$
with $p=p_N$. Then $e_N\leq p_N$, contradicting
$E_N\cap\downarrow_F p_N=\varnothing$.

Hence $F$ has a greatest element $p$ and then
$$F=\downarrow_P p=\operatorname{cl}_{\Sigma P}\{p\}.$$ Therefore $\Sigma P$ is sober.
\end{proof}

\begin{theorem}
Let $P$ be a countable dcpo satisfying \textup{(FG)}. If $\Sigma P$
is Fr\'echet--Urysohn, then $\Sigma P$ is sober.
\end{theorem}

\begin{proof}
By Proposition~3.3, $P$ has weak one-step closure. Since $P$ is
countable, Theorem~\ref{thm:countable-weak-one-step-sober} gives
sobriety of $\Sigma P$.
\end{proof}
\section{A countable non-sober Fr\'echet--Urysohn Scott space}\label{x:section}

We construct a countable dcpo whose Scott space is Fr\'echet--Urysohn but not sober.

Write $\mathbb N_0=\{0,1,\ldots\}$ and $\mathbb N_+=\{1,2,\ldots\}$.
Let $P$ be the union of the following three families of distinct symbols:
$\{m_j^k:k\in\mathbb N_0,\ j\in\mathbb N_+\}$,
$\{d_{j,l}^k:k\in\mathbb N_0,\ j,l\in\mathbb N_+\}$, and
$\{a_r^k:k,r\in\mathbb N_+\}$. In particular, there is no
$a$-element at layer zero. Let $\leq$ be the reflexive and transitive
closure of the following relations.

\smallskip
\noindent\textup{(R1)}\quad
$d_{j,l}^k\leq d_{j,l+1}^k\leq m_j^k$
for $k\in\mathbb N_0$ and $j,l\in\mathbb N_+$.

\noindent\textup{(R2)}\quad
$d_{j,l}^k\leq m_r^{k+1}$
for $k\in\mathbb N_0$ and $j,l,r\in\mathbb N_+$ with $l\leq r$.

\noindent\textup{(R3)}\quad
$d_{j,l}^k\leq a_r^k$
for $k,j,l,r\in\mathbb N_+$ with $l\leq r$.

\noindent\textup{(R4)}\quad
$a_r^k\leq d_{j,l}^{k-1}$
for $k,r,j,l\in\mathbb N_+$ with $r\leq j$ and $r\leq l$.
\smallskip

For $k\in\mathbb N_0$ and $j\in\mathbb N_+$, put
$D_j^k=\{d_{j,l}^k:l\in\mathbb N_+\}$. The superscript of a
$d$- or $a$-element is its \emph{layer}. Define its \emph{height} by
$w(d_{j,l}^k)=l$ and $w(a_r^k)=r$.

Consider a path $x_0,\ldots,x_n$ of nonidentity generating steps
containing no $m$-element, and let $k_i$ be the layer of $x_i$.
A step within a row $D_j^k$ preserves the layer and increases the
height from $l$ to $l+1$. An \textup{(R3)} step preserves the layer
and changes the height from $l$ to $r\geq l$. An \textup{(R4)} step
changes the layer from $k$ to $k-1$ and the height from $r$ to
$l\geq r$. Hence $k_{i+1}\leq k_i$ and
$w(x_i)\leq w(x_{i+1})$ for every $i<n$.

Suppose that $n\geq1$ and $x_n=x_0$. Then
$k_0\geq\cdots\geq k_n=k_0$ and
$w(x_0)\leq\cdots\leq w(x_n)=w(x_0)$.
Thus the path contains neither an \textup{(R4)} step nor a step
within a row. The only remaining possibility is an \textup{(R3)}
step ending at an $a$-element. Its next step around the cycle would
have to be an \textup{(R4)} step, a contradiction. A nonidentity
cycle containing an $m$-element is also impossible, since no
nonidentity generating step starts at an $m$-element.

If $x\leq y\leq x$ with $x\neq y$, concatenating a path from $x$ to
$y$ with a path from $y$ to $x$ gives a nonidentity cycle. Hence
$\leq$ is antisymmetric and therefore a partial order on $P$.
Moreover, $\uparrow m_j^k=\{m_j^k\}$ for every $k\in\mathbb N_0$
and $j\in\mathbb N_+$, so each $m_j^k$ is maximal. Finally, each
of the three defining families is indexed by a finite product of
countable sets, and $\{m_j^0:j\in\mathbb N_+\}$ is infinite.
Thus $P$ is countably infinite.

\begin{lemma}\label{lem:layer-height-bound}
Suppose that $d_{s,l}^h\leq m_j^k$ and $(h,s)\neq(k,j)$.
Then $k\leq h+1$ and $l\leq j$. If $d_{s,l}^h\leq z$ for a
nonmaximal element $z$, then $l\leq w(z)$.
\end{lemma}

\begin{proof}
No nonidentity generating step starts at an $m$-element. Hence a
path from $d_{s,l}^h$ to a nonmaximal element $z$ contains only
$d$- and $a$-elements. Along such a path the layer is nonincreasing
and the height is nondecreasing. In particular,
$l=w(d_{s,l}^h)\leq w(z)$.

Now take a path from $d_{s,l}^h$ to $m_j^k$ and consider its last
nonidentity step. If it is an \textup{(R2)} step, then $k\geq1$ and
the step has the form $d_{i,t}^{k-1}\leq m_j^k$ with $t\leq j$.
The preceding path gives $k-1\leq h$ and $l\leq t$. Thus
$k\leq h+1$ and $l\leq t\leq j$.

If the last step is an \textup{(R1)} step, its initial point belongs
to $D_j^k$. Since $(h,s)\neq(k,j)$, the path starts outside $D_j^k$.
Its first entry into $D_j^k$ must be an \textup{(R4)} step:
an \textup{(R1)} step ending in $D_j^k$ starts in that same row.
Thus the path contains a step $a_r^{k+1}\leq d_{j,t}^k$ with
$r\leq j$ and $r\leq t$. The path from $d_{s,l}^h$ to $a_r^{k+1}$
gives $k+1\leq h$ and $l\leq r$. In particular, $k\leq h+1$ and
$l\leq r\leq j$.
\end{proof}

\begin{proposition}\label{prop:layered-dcpo}
The poset $P$ is a dcpo. Every directed subset $D$ of $P$ has a
greatest element or has a cofinal subset of unbounded height in
some row $D_j^k$. In the latter case $\bigvee D=m_j^k$.
In particular, $\bigvee D_j^k=m_j^k$ for every $k\in\mathbb N_0$
and $j\in\mathbb N_+$.
\end{proposition}

\begin{proof}
Let $D\subseteq P$ be directed. If $m_j^k\in D$, then for every
$z\in D$ there exists $u\in D$ with $z,m_j^k\leq u$. Maximality
gives $u=m_j^k$, so $m_j^k$ is the greatest element of $D$.

Suppose that $D$ contains no $m$-element, and let $h$ be the least
layer represented in $D$. If $a_r^h\in D$, then every nonidentity
path from $a_r^h$ to a nonmaximal element begins with an
\textup{(R4)} step to layer $h-1$. Its subsequent layers are at most
$h-1$. Hence $D\cap\uparrow a_r^h=\{a_r^h\}$. For each $z\in D$,
a common upper bound of $z$ and $a_r^h$ in $D$ must therefore be
$a_r^h$. Thus $a_r^h$ is the greatest element of $D$.

It remains to consider the case in which all points of $D$ at layer
$h$ are row elements. Choose $d_{j,l_0}^h\in D$ and put
$D'=D\cap\uparrow d_{j,l_0}^h$. By directedness, each $z\in D$
has an upper bound in $D'$. Any common upper bound in $D$ of two
members of $D'$ also belongs to $D'$. Thus $D'$ is directed and
cofinal in $D$.

For $z\in D'$, a path from $d_{j,l_0}^h$ to $z$ cannot increase
the layer. By minimality of $h$, the path stays at layer $h$.
If it uses an \textup{(R3)} step, it ends at an $a$-element or
leaves layer $h$ by \textup{(R4)}. Both possibilities are excluded.
Hence $D'\subseteq D_j^h$. If its heights are bounded, the element
of $D'$ of largest height is greatest in $D'$ and, by cofinality,
in $D$.

Suppose that the heights in $D'$ are unbounded. By \textup{(R1)},
$m_j^h$ is an upper bound of $D'$. If $z$ were a nonmaximal upper
bound, Lemma~\ref{lem:layer-height-bound} would give $l\leq w(z)$
for every $d_{j,l}^h\in D'$. If $m_i^k\neq m_j^h$ were an upper
bound, the same lemma would give $l\leq i$ for every
$d_{j,l}^h\in D'$. Both conclusions contradict unboundedness.
Every maximal element is an $m$-element, since
$d_{j,l}^k<d_{j,l+1}^k$ and $a_r^k<d_{r,r}^{k-1}$.
Thus $m_j^h$ is the unique upper bound of $D'$. Since $D'$ is
cofinal in $D$, $\bigvee D=m_j^h$. Taking $D=D_j^k$ gives
$\bigvee D_j^k=m_j^k$.
\end{proof}

\begin{lemma}\label{lem:layer-tails}
Every nonempty Scott-open set contains a tail of $(m_j^k)_{j\geq1}$
for each $k\in\mathbb N_0$ and a tail of $(a_r^k)_{r\geq1}$ for
each $k\in\mathbb N_+$. In particular, $a_r^k\to x$ as $r\to\infty$
for every $k\in\mathbb N_+$ and every $x\in P$.
\end{lemma}

\begin{proof}
Let $\varnothing\neq U\in\sigma(P)$. Since $d_{j,l}^k\leq m_j^k$
and $a_r^k\leq d_{r,r}^{k-1}\leq m_r^{k-1}$, upperness gives
$m_j^k\in U$ for some $k,j$. Whenever $m_j^k\in U$,
Proposition~\ref{prop:layered-dcpo} gives $m_j^k=\bigvee D_j^k$.
Hence $d_{j,l}^k\in U$ for some $l\geq1$.

By \textup{(R2)}, $d_{j,l}^k\leq m_r^{k+1}$ for every $r\geq l$,
so $m_r^{k+1}\in U$ for every $r\geq l$. If $k\geq1$, then
\textup{(R3)} also gives $d_{j,l}^k\leq a_r^k$ for every $r\geq l$.
In particular, $a_l^k\in U$. For every $i\geq l$,
\textup{(R4)} and \textup{(R1)} give
$a_l^k\leq d_{i,l}^{k-1}\leq m_i^{k-1}$, whence $m_i^{k-1}\in U$.

Thus a maximal element in $U$ at layer $k$ yields a tail of maximal
elements at layer $k+1$ and, when $k\geq1$, at layer $k-1$.
Iterating downwards to layer zero and then upwards shows that $U$
contains a maximal element at every layer. Applying the upward step
from layer $k-1$ gives a tail of $(m_j^k)_{j\geq1}$ for each $k\geq1$.
The downward step from layer one gives the tail at layer zero.
At every positive layer $k$, the \textup{(R3)} step above gives a
tail of $(a_r^k)_{r\geq1}$. This applies to every Scott-open
neighbourhood of every $x\in P$, proving $a_r^k\to x$.
\end{proof}

\begin{proposition}\label{prop:layered-nonsober}
The Scott space $\Sigma P$ is irreducible and non-sober.
\end{proposition}

\begin{proof}
Let $U,V\in\sigma(P)$ be nonempty. By Lemma~\ref{lem:layer-tails},
there exist $r_U,r_V\geq1$ such that $a_r^1\in U$ for $r\geq r_U$
and $a_r^1\in V$ for $r\geq r_V$. Hence
$a_{\max\{r_U,r_V\}}^1\in U\cap V$. Thus $P$ is an irreducible
Scott-closed set.

If $p$ were a generic point of $P$, then
$P=\operatorname{cl}_{\Sigma P}\{p\}=\downarrow p$. In particular,
$m_1^0\leq p$ and $m_2^0\leq p$. Maximality would give
$m_1^0=p=m_2^0$, contrary to $m_1^0\neq m_2^0$. Hence $P$ has no
generic point, and $\Sigma P$ is non-sober.
\end{proof}

\begin{proposition}\label{prop:layered-closure}
Let $L\subseteq P$ be lower and, for $k\in\mathbb N_0$, put
$J_k=\{j\in\mathbb N_+:D_j^k\subseteq L\}$. If some $J_k$ is
unbounded, then $\operatorname{cl}_{\Sigma P}L=P$. If every $J_k$
is bounded, then
$\operatorname{cl}_{\Sigma P}L
 =L\cup\bigcup_{k\in\mathbb N_0}\bigcup_{j\in J_k}\downarrow m_j^k$.
\end{proposition}

\begin{proof}
Suppose that $J_k$ is unbounded. For each $r\geq1$, choose
$j\in J_k$ with $j\geq r$. By \textup{(R4)},
$a_r^{k+1}\leq d_{j,r}^k\in D_j^k\subseteq L$. Since $L$ is lower,
$a_r^{k+1}\in L$ for every $r\geq1$. By Lemma~\ref{lem:layer-tails},
$a_r^{k+1}\to x$ for every $x\in P$. Hence
$\operatorname{cl}_{\Sigma P}L=P$.

Suppose that every $J_k$ is bounded, hence finite. Put
$H=L\cup\bigcup_{k\in\mathbb N_0}\bigcup_{j\in J_k}\downarrow m_j^k$.
For $j\in J_k$, we have $D_j^k\subseteq L$ and
$m_j^k=\bigvee D_j^k\in\operatorname{cl}_{\Sigma P}L$.
Thus $H$ is lower and $L\subseteq H\subseteq\operatorname{cl}_{\Sigma P}L$.

Fix $h\in\mathbb N_0$ and $s\notin J_h$. Choose $t\geq1$ with
$d_{s,t}^h\notin L$. Since $L$ is lower, $d_{s,l}^h\in L$ implies
$l<t$. If $d_{s,l}^h\in H\setminus L$, then
$d_{s,l}^h\leq m_j^k$ for some $k\in\mathbb N_0$ and $j\in J_k$.
Since $s\notin J_h$, we have $(k,j)\neq(h,s)$.
Lemma~\ref{lem:layer-height-bound} gives $k\leq h+1$ and $l\leq j$.
Consequently, every $d_{s,l}^h\in H$ satisfies $l\leq b$, where
$b=\max\bigl(\{t-1\}\cup\bigcup_{k=0}^{h+1}J_k\bigr)$.
This maximum is finite, so $d_{s,b+1}^h\notin H$. We have proved
that $D_s^h\subseteq H$ if and only if $s\in J_h$.

Let $D\subseteq H$ be directed. If $D$ has a greatest element, then
$\bigvee D\in D\subseteq H$. Otherwise,
Proposition~\ref{prop:layered-dcpo} gives a cofinal subset of $D$
of unbounded height in some row $D_s^h$, with $\bigvee D=m_s^h$.
For every $l\geq1$, this subset contains some $d_{s,t}^h$ with
$t\geq l$. Thus $d_{s,l}^h\leq d_{s,t}^h\in H$, and lowerness gives
$D_s^h\subseteq H$. Hence $s\in J_h$ and $\bigvee D=m_s^h\in H$.
Therefore $H$ is Scott closed. Since
$L\subseteq H\subseteq\operatorname{cl}_{\Sigma P}L$, it follows
that $\operatorname{cl}_{\Sigma P}L=H$.
\end{proof}

\begin{theorem}\label{thm:layered-FU}
The Scott space $\Sigma P$ is Fr\'echet--Urysohn.
\end{theorem}

\begin{proof}
Let $A\subseteq P$ and $x\in\operatorname{cl}_{\Sigma P}A$.
Put $L=\downarrow A$. Since $A\subseteq L\subseteq
\operatorname{cl}_{\Sigma P}A$, we have
$\operatorname{cl}_{\Sigma P}L=\operatorname{cl}_{\Sigma P}A$.
Define $J_k=\{j\in\mathbb N_+:D_j^k\subseteq L\}$.

Suppose first that some $J_k$ is unbounded. As in the proof of
Proposition~\ref{prop:layered-closure}, $a_r^{k+1}\in L$ for every
$r\geq1$. Choose $b_r\in A$ with $a_r^{k+1}\leq b_r$.
For each $x\in U\in\sigma(P)$, Lemma~\ref{lem:layer-tails} gives
$R\geq1$ such that $a_r^{k+1}\in U$ for every $r\geq R$.
By upperness, $b_r\in U$ for every $r\geq R$. Thus $b_r\to x$.

Suppose next that every $J_k$ is bounded. If $x\in L$, choose
$b\in A$ with $x\leq b$. Every Scott-open neighbourhood of $x$
contains $b$, so the constant sequence at $b$ converges to $x$.

Finally, let $x\notin L$. By Proposition~\ref{prop:layered-closure},
there exist $k\in\mathbb N_0$ and $j\in J_k$ with $x\leq m_j^k$.
Since $D_j^k\subseteq L=\downarrow A$, choose $b_l\in A$ with
$d_{j,l}^k\leq b_l$ for every $l\geq1$. If $x\in U\in\sigma(P)$,
then $m_j^k\in U$. Since $m_j^k=\bigvee D_j^k$, there exists $l_0$
with $d_{j,l_0}^k\in U$. For every $l\geq l_0$ we have
$d_{j,l_0}^k\leq d_{j,l}^k\leq b_l$, whence $b_l\in U$.
Thus $b_l\to x$, as required.
\end{proof}

\section{Every meet-continuous locally compact dcpo is Scott sober}\label{sec:local}

The preceding results use the countability of the underlying order. We
now consider a sufficient condition involving the order and its Scott
topology that applies to dcpos of arbitrary cardinality.

\begin{lemma}\label{lem:web-open} {\rm (\cite{domains}) Let $P$ be a dcpo. Then $P$ is meet continuous if and only if for each $U\in\sigma(P)$ and $x\in U$, $\up(U\cap\dn x)\in\sigma(P)$.}
\end{lemma}

For a dcpo $X$ and a subset $A\subseteq X$, put
\begin{equation}\label{eq:h}
                         h_X(A)=\Int_{\Sigma X}(\up_X A).
\end{equation}

\begin{lemma}
{\rm Let $P$ be a meet-continuous dcpo. For $W\in\sigma(P)$ and
$r\geq 1$, the set
\begin{equation}\label{eq:gamma-r}
\Gamma_r(W)
=
\bigl\{(y_1,\ldots,y_r)\in P^r:
       \exists w\in W,\forall i\in\{1,\ldots,r\},\ w\leq y_i
\bigr\}
\end{equation}
is Scott open in $P^r$.}
\end{lemma}

\begin{proof}
If $(y_1,\ldots,y_r)\in\Gamma_r(W)$ and
$(y_1,\ldots,y_r)\leq(z_1,\ldots,z_r)$, we can choose $w\in W$ with
$w\leq y_i$ for every $i$. Then $w\leq y_i\leq z_i$ for every $i$.
So $(z_1,\ldots,z_r)\in\Gamma_r(W)$. Thus $\Gamma_r(W)$ is upper.

Let $\mathcal D\subseteq P^r$ be a directed subset with
$\bigvee\mathcal D\in\Gamma_r(W)$. For $1\leq i\leq r$, let
$\pi_i:P^r\to P$ be the $i$th projection and put
$t_i=\bigvee\pi_i[\mathcal D]$. Since directed suprema in $P^r$
are coordinatewise,
\[
\bigvee\mathcal D=(t_1,\ldots,t_r)\in\Gamma_r(W).
\]
Choose $w_0\in W$ such that $w_0\leq t_i$ for every $i$.

Inductively, suppose that $1\leq i\leq r$ and that
$w_{i-1}\in W$ has been chosen with $w_{i-1}\leq w_0$.
By Lemma~5.1, $\uparrow(W\cap\downarrow w_{i-1})\in\sigma(P)$. Since $w_{i-1}\in W\cap\downarrow w_{i-1}$ and
$w_{i-1}\leq w_0\leq t_i$, we have
\[
\bigvee\pi_i[\mathcal D]=t_i
   \in\uparrow(W\cap\downarrow w_{i-1}).
\]
The set $\pi_i[\mathcal D]$ is directed. Hence there exists
$\delta^i\in\mathcal D$ such that $\pi_i(\delta^i)\in\uparrow(W\cap\downarrow w_{i-1})$. Choose $w_i\in W\cap\downarrow w_{i-1}$ with
$w_i\leq\pi_i(\delta^i)$. This gives
\[
w_r\leq\cdots\leq w_1\leq w_0,\
w_i\in W\mbox{ and } w_i\leq\pi_i(\delta^i)
(1\leq i\leq r).
\]

By directedness, choose $\delta\in\mathcal D$ with
$\delta^i\leq\delta$ for all $1\leq i\leq r$. Then $w_r\in W$ and
$w_r\leq w_i\leq\pi_i(\delta^i)\leq\pi_i(\delta)(1\leq i\leq r)$. Thus $\delta\in\mathcal D\cap\Gamma_r(W)$, proving that
$\Gamma_r(W)\in\sigma(P^r)$.
\end{proof}

\begin{lemma}\label{lem:product}{\rm(\cite{domains})
If $S$ is a dcpo and $\Sigma P$ is core-compact, then $\Sigma (S\times P)=\Sigma S\times\Sigma P$.}
\end{lemma}

\begin{lemma}\label{lem:local-arity}
{\rm Let $P$ be a locally compact meet-continuous dcpo. If
$y\in W\in\sigma(P)$ and $r\geq1$, then there exists
$T\in\sigma(P)$ with $y\in T\subseteq W$ such that every nonempty
subset of $T$ of size at most $r$ has a common lower bound in $W$.}
\end{lemma}

\begin{proof}
As $P$ is a locally compact dcpo, $\Sigma P$ is core-compact. For every $n\geq1$, by Lemma~5.3,
\[
\Sigma(P^{n+1})=\Sigma(P^n)\times\Sigma P.
\]
Hence, $\Sigma(P^r)=(\Sigma P)^r$. By Lemma~5.2,
$\Gamma_r(W)$ is open in $(\Sigma P)^r$. Since $y\in W$,
$(y,\ldots,y)\in\Gamma_r(W)$.

Choose $T_1,\ldots,T_r\in\sigma(P)$ such that
\[
(y,\ldots,y)\in T_1\times\cdots\times T_r
                 \subseteq\Gamma_r(W).
\]
Put $T=W\cap\bigcap_{i=1}^r T_i$. Then
$T\in\sigma(P)$, $y\in T\subseteq W$ and
$T^r\subseteq T_1\times\cdots\times T_r\subseteq\Gamma_r(W)$.

Let $F=\{x_1,\ldots,x_m\}\subseteq T$, where $1\leq m=|F|\leq r$.
Then
\[
\bigl(x_1,\ldots,x_m,
      \underbrace{x_m,\ldots,x_m}_{r-m}\bigr)
\in T^r\subseteq\Gamma_r(W).
\]
By the definition of $\Gamma_r(W)$, there exists $w\in W$ such that
$w\leq x_i$ for every $1\leq i\leq m$. Thus $w$ is a common lower
bound of $F$ in $W$.
\end{proof}

\begin{lemma}\label{lem:avoid}
{\rm Let $P$ be a locally compact meet-continuous dcpo and write $h=h_P$.
If $E\subseteq P$ and $x\in h(E)$, then there is a Scott open
neighbourhood $V$ of $x$ such that $V\subseteq\uparrow(E\cap\downarrow v)$.
If $N\subseteq P$ and $x\notin h(N)$, then there exists
a relatively closed subset $E'$ of $E$ such that
$E'\cap N=\varnothing$ and $x\in h(E')$.}
\end{lemma}

\begin{proof}
Since $x\in h(E)=\operatorname{int}_{\Sigma P}(\up E)$,
by Lemma~5.4, there is $V\in\sigma(P)$
with $x\in V\subseteq h(E)$ such that
\[
\forall u,v\in V,\exists w\in h(E),
 w\in\dn u\cap\dn v.
\]
Fix $v\in V$. For each $u\in V$, choose such a $w$. Since
$w\in h(E)\subseteq\uparrow E$, there exists $e\in E$ with
$e\leq w$. Thus $e\leq w\leq v$ and $e\leq w\leq u$.
Consequently, for every $v\in V$, $x\in V\subseteq h(E\cap\dn v)$ if $V\subseteq\up(E\cap\dn v)$.
Now suppose that $x\notin h(N)$. If $V\subseteq\uparrow N$, then
$x\in V\subseteq\operatorname{int}_{\Sigma P}(\uparrow N)=h(N)$.
Hence there exists $v\in V\setminus\uparrow N$. Put
\[
E'=E\cap\downarrow v.
\]
Then $E'$ is relatively
closed in $E$. As $v\notin\uparrow N$, we know $E'\cap N=\varnothing$. Thus $x\in V\subseteq h(E')$.
\end{proof}

\label{lem:finite-union}

\begin{lemma}
Let $P$ be a locally compact meet-continuous dcpo and write $h=h_P$.
Then for every family $(A_k)_{1\leq k\leq n}$ of subsets of $P$,
\[
h\Bigl(\bigcup_{1\leq k\leq n}A_k\Bigr)=\bigcup_{1\leq k\leq n}h(A_k).
\]
\end{lemma}

\begin{proof}
For $A\subseteq B\subseteq P$, we have
$h(A)=\operatorname{int}_{\Sigma P}(\up A)
\subseteq\operatorname{int}_{\Sigma P}(\up B)=h(B)$. Thus, for all $A,B\subseteq P$,
$h(A)\cup h(B)\subseteq h(A\cup B)$.

Conversely, let $x\in h(A\cup B)\setminus h(A)$.
By Lemma~5.5, there exists a
relatively closed subset $E'$ of $A\cup B$ such that
\[
E'\cap A=\varnothing
\text{ and }
x\in h(E').
\]
Since $E'\subseteq(A\cup B)\setminus A\subseteq B$, monotonicity gives
$x\in h(E')\subseteq h(B)$. Hence
$h(A\cup B)\setminus h(A)\subseteq h(B)$ and therefore
$h(A\cup B)=h(A)\cup h(B)$. It is not difficult to verify that the statement also holds for $n$ sets.
\end{proof}

\begin{lemma}\label{lem:closed-refinement}
{\rm Let $P$ be a locally compact meet-continuous dcpo and write $h=h_P$.
Let $D\subseteq P$ be a compact subset and $\mathcal U$
 a relatively open cover of $D$. Then
\begin{equation}\label{eq:closed-refinement}
h(D)=
\bigcup\bigl\{h(C): C\text{ is closed in }D,\
                    C\subseteq U\text{ for some }U\in\mathcal U\bigr\}.
\end{equation}}
\end{lemma}

\begin{proof}
For each $C\subseteq D$, monotonicity gives $h(C)\subseteq h(D)$.
Thus the right-hand side of (5.3) is contained
in $h(D)$.

Fix $x\in h(D)$. By compactness, choose $U_1,\ldots,U_m\in\mathcal U$
such that $D=\bigcup_{j=1}^m U_j$. Put
\[
A_j=D\setminus U_j,\
B_0=D\text{ and }
B_j=B_{j-1}\cap A_j
   =D\setminus\bigcup_{i=1}^j U_i(1\leq j\leq m).
\]
Each $B_j$ is closed in $D$. Since $x\in h(B_0)$ and
$h(B_m)=h(\varnothing)=\varnothing$, the index
\[
j=\min\{i\in\{1,\ldots,m\}:x\notin h(B_i)\}
\]
is well defined and satisfies $x\in h(B_{j-1})\setminus h(B_j)$. By Lemma~5.5, there exists
$C\subseteq B_{j-1}$ closed in $B_{j-1}$ such that
\[
x\in h(C), C\cap B_j=\varnothing.
\]
Since $B_{j-1}$ is closed in $D$, so is $C$. Moreover,
$C\subseteq B_{j-1}\setminus B_j
 =B_{j-1}\cap U_j\subseteq U_j$. Thus $h(C)$ occurs in the union in (5.3) and
contains $x$. This proves the formula (5.3).
\end{proof}

\begin{lemma}\label{lem:prefix}
{\rm Let $P$ be a locally compact meet-continuous dcpo and write $h=h_P$.
Let $r\geq1$, $m\geq0$, and let $L_0,\ldots,L_m$ be nonempty
compact subsets of $P$ such that for each $0\leq i<m$, $L_{i+1}\subseteq h(L_i)$. For $0\leq j\leq m$, let $\mathcal D_j$ be the family of all lists
$(D_0,\ldots,D_j)$ satisfying the following conditions.
\begin{enumerate}
\item $D_0=L_0$ and $D_i$ is a nonempty closed subset of $L_i$
      for $0\leq i\leq j$.
\item $D_{i+1}\subseteq h(D_i)$ for $0\leq i<j$.
\item For $0\leq i<j$ and $F\subseteq D_{i+1}$ with
      $1\leq |F|\leq r$, $D_i\cap\bigcap_{v\in F}\downarrow v\neq\varnothing$.

\end{enumerate}
Then, for every $0\leq j\leq m$,
\begin{equation}\label{eq:compact-lists}
h(L_j)=
\bigcup_{(D_0,\ldots,D_j)\in\mathcal D_j}h(D_j).
\end{equation}
In particular, for each $z\in h(L_m)$ there exists
$(D_0,\ldots,D_m)\in\mathcal D_m$ such that $z\in h(D_m)$.}
\end{lemma}

\begin{proof}
For each $(D_0,\ldots,D_j)\in\mathcal D_j$, monotonicity gives
$h(D_j)\subseteq h(L_j)$. It remains to prove the reverse inclusion.
For $j=0$, we have $\mathcal D_0=\{(L_0)\}$. So
\eqref{eq:compact-lists} holds.

Suppose that \eqref{eq:compact-lists} holds for some $j<m$. Then
\[
L_{j+1}\subseteq h(L_j)
=\bigcup_{(D_0,\ldots,D_j)\in\mathcal D_j}h(D_j).
\]
For each $y\in L_{j+1}$, choose
$(D_0^y,\ldots,D_j^y)\in\mathcal D_j$ with $y\in h(D_j^y)$.
By Lemma~5.4, there exists $T_y\in\sigma(P)$ such that
$y\in T_y\subseteq h(D_j^y)$ and for $F\subseteq T_y,\ 1\leq |F|\leq r$,
$h(D_j^y)\cap\bigcap\limits_{v\in F}\dn v\neq\varnothing$.

For each such $F$, choose $w\in h(D_j^y)$ with $w\leq v$ for all
$v\in F$. Since $h(D_j^y)\subseteq\uparrow D_j^y$, there exists
$d\in D_j^y$ with $d\leq w$. Thus
$d\leq w\leq v\quad(v\in F)$ and further $D_j^y\cap\bigcap_{v\in F}\downarrow v\neq\varnothing$.
Now fix $z\in h(L_{j+1})$. The family
\[
\mathcal U=\{L_{j+1}\cap T_y:y\in L_{j+1}\}
\]
is a relatively open cover of $L_{j+1}$. By Lemma~5.7, there exist
$y\in L_{j+1}$ and a closed subset $D_{j+1}$ of $L_{j+1}$ such that
\[
z\in h(D_{j+1})\mbox{ and }
D_{j+1}\subseteq L_{j+1}\cap T_y\subseteq h(D_j^y).
\]
Since $h(\varnothing)=\varnothing$, we have $D_{j+1}\neq\varnothing$.
Moreover, $D_{j+1}\subseteq T_y$ gives
\[
D_j^y\cap\bigcap_{v\in F}\dn v\neq\varnothing (F\subseteq D_{j+1},\ 1\leq |F|\leq r).
\]
Hence $(D_0^y,\ldots,D_j^y,D_{j+1})\in\mathcal D_{j+1}$ and
$z\in h(D_{j+1})$. This proves the reverse inclusion at $j+1$
and completes the induction. The last assertion follows from
\eqref{eq:compact-lists} with $j=m$.
\end{proof}

\begin{lemma}\label{lem:diagonal-compactness}
{\rm Let $P$ be a locally compact meet-continuous dcpo. Suppose $K_i$ is a nonempty compact subset of $P$ for each $i\geq0$.
For $i\leq k$, let $F_i^{(k)}$ be a nonempty closed subset of
$K_i$ such that $F_i^{(k+1)}\subseteq F_i^{(k)}$. Suppose for
$0\leq i<k$, $$F_i^{(k)}\cap\bigcap\limits_{y\in B}\dn y\neq\varnothing$$ if $B\subseteq F_{i+1}^{(k)}$ and $1\leq |B|\leq k$. Put $E_i=\bigcap\limits_{k\geq i}F_i^{(k)}$. Then each $E_i$ is nonempty and one can choose $p_i\in E_i$ such that $p_i\leq y$ for every
$y\in E_{i+1}$. Furthermore, any such choices satisfy
$p_i\in K_i\mbox{ and } p_i\leq p_{i+1}(i\geq0)$.
}
\end{lemma}

\begin{proof}
Fix $i\geq0$. For every nonempty finite set
$J\subseteq\{i,i+1,\ldots\}$, we have
\[
\bigcap_{j\in J}F_i^{(j)}=F_i^{(\max J)}\neq\varnothing.
\]
Since $K_i$ is compact and each $F_i^{(j)}$ is closed in $K_i$,
$E_i=\bigcap_{k\geq i}F_i^{(k)}\neq\varnothing$. Clearly, the members of
\[
\bigl\{F_i^{(j)}:j\geq i\bigr\}
\ \cup\
\bigl\{K_i\cap\downarrow y:y\in E_{i+1}\bigr\}
\]
are closed in $K_i$.
Take finite sets $J\subseteq\{i,i+1,\ldots\}$ and $B\subseteq E_{i+1}$.
Put $k=\max\bigl(J\cup\{i+1,|B|\}\bigr)$. Then $i<k$, $|B|\leq k$ and
\[
B\subseteq E_{i+1}\subseteq F_{i+1}^{(k)},
F_i^{(k)}\subseteq F_i^{(j)}\quad(j\in J).
\]
If $B\neq\varnothing$, then $a\in F_i^{(k)}$ with
$a\leq y$ for every $y\in B$. If $B=\varnothing$, take any
$a\in F_i^{(k)}$. In either case,
\[
a\in K_i,
a\in F_i^{(j)}(j\in J),
a\in K_i\cap\downarrow y(y\in B).
\]
Compactness of $K_i$ therefore gives
\[
\varnothing\neq
\left(\bigcap_{j\geq i}F_i^{(j)}\right)
\cap\bigcap_{y\in E_{i+1}}(K_i\cap\downarrow y)
=
E_i\cap\bigcap_{y\in E_{i+1}}\downarrow y.
\]
For each $i\geq0$, choose $p_i$ in this intersection. Then
$p_i\in E_i\subseteq K_i$ and $p_i\leq y$ for every $y\in E_{i+1}$.
Taking $y=p_{i+1}\in E_{i+1}$ gives $p_i\leq p_{i+1}$.
\end{proof}

\begin{theorem}\label{thm:selection}
{\rm Let $P$ be a meet-continuous locally compact dcpo. Suppose that $(K_n)_{n\in\mathbb N_0}$ is a sequence of nonempty compact subsets of $P$ satisfying $K_{n+1}\subseteq h(K_n)(n\in\mathbb N_0)$.
Then there exist points $p_n\in K_n$ such that
$p_n\leq p_{n+1}$ for every $n\in\mathbb N_0$.}
\end{theorem}

\begin{proof}
We construct a rooted tree with a finite nonempty set $\mathcal T_k$
of nodes at each level $k$. Each $v\in\mathcal T_k$ carries a list
$(F_0(v),\ldots,F_k(v))$, where $F_i(v)$ is a nonempty closed subset
of $K_i$. We require

(1) $F_{i+1}(v)\subseteq h(F_i(v))$ for $v\in\mathcal T_k,\ 0\leq i<k$;

(2) $F_i(v)\cap\bigcap_{b\in B}\downarrow b\neq\varnothing$ for $v\in\mathcal T_k,\ 0\leq i<k,\
      B\subseteq F_{i+1}(v),\ 1\leq |B|\leq k$;

(3) $K_{k+1}\subseteq\bigcup_{v\in\mathcal T_k}h(F_k(v))$.

If $v'\in\mathcal T_{k+1}$ is a child of $v\in\mathcal T_k$, we also
require $F_i(v')\subseteq F_i(v)$ for $0\leq i\leq k$.
At level $0$, take a single root carrying $(K_0)$. The covering
condition holds since $K_1\subseteq h(K_0)$.

Suppose that $\mathcal T_k$ has been constructed. For each
$v\in\mathcal T_k$, using Lemma~5.8 with
\[
m=k,\ r=k+1,\ L_i=F_i(v)(0\leq i\leq k).
\]
Indeed, each $F_i(v)$ is closed in the compact space $K_i$ and hence
is compact. Writing $\mathcal D_k(v)$ for the family of lists supplied
by that lemma, we have
\[
K_{k+1}
\subseteq\bigcup_{v\in\mathcal T_k}h(F_k(v))
=\bigcup_{v\in\mathcal T_k}
  \ \bigcup_{(D_0,\ldots,D_k)\in\mathcal D_k(v)}h(D_k).
\]
For each $y\in K_{k+1}$, choose $v(y)\in\mathcal T_k$ and
\[
(D_0^y,\ldots,D_k^y)\in\mathcal D_k(v(y))
\text{ with } y\in h(D_k^y).
\]
By Lemma~5.4, there exists $T_y\in\sigma(P)$ such that
$y\in T_y\subseteq h(D_k^y)$ and
\[
h(D_k^y)\cap\bigcap_{b\in B}\downarrow b\neq\varnothing (B\subseteq T_y,\ 1\leq |B|\leq k+1).
\]
For any such $B$, choose $w$ in this intersection. Since
$w\in h(D_k^y)\subseteq\uparrow D_k^y$, some $d\in D_k^y$ satisfies
$d\leq w\leq b$ for every $b\in B$. Hence
\[
D_k^y\cap\bigcap_{b\in B}\downarrow b\neq\varnothing(B\subseteq T_y,\ 1\leq |B|\leq k+1).
\]

The family $\{K_{k+1}\cap T_y:y\in K_{k+1}\}$ is a relatively open
cover of $K_{k+1}$. By Lemma~5.7,
\[
h(K_{k+1})=
\bigcup\bigl\{h(C):
   C\text{ is closed in }K_{k+1},\
   C\subseteq T_y\text{ for some }y\in K_{k+1}\bigr\}.
\]
Since $K_{k+2}\subseteq h(K_{k+1})$ is nonempty and compact and
$h(\varnothing)=\varnothing$, there exist $a\geq1$,
$y_1,\ldots,y_a\in K_{k+1}$ and nonempty closed subsets
$C_1,\ldots,C_a$ of $K_{k+1}$ such that
\[
C_\ell\subseteq T_{y_\ell} (1\leq\ell\leq a)\mbox{ and }
K_{k+2}\subseteq\bigcup_{\ell=1}^a h(C_\ell).
\]
For each $\ell$, adjoin a child of $v(y_\ell)$ carrying the list
$(D_0^{y_\ell},\ldots,D_k^{y_\ell},C_\ell)$. For $0\leq i\leq k$, the set $D_i^{y_\ell}$ is closed in
$F_i(v(y_\ell))$ and hence closed in $K_i$. For $i<k$, the first two
requirements at level $k+1$ follow from
$(D_0^{y_\ell},\ldots,D_k^{y_\ell})\in\mathcal D_k(v(y_\ell))$
with $r=k+1$. For $i=k$, they follow from
\[
C_\ell\subseteq T_{y_\ell}\subseteq h(D_k^{y_\ell})\mbox{ and }
D_k^{y_\ell}\cap\bigcap_{b\in B}\downarrow b\neq\varnothing
(B\subseteq C_\ell,\ 1\leq |B|\leq k+1).
\]
The sets $h(C_\ell)$ cover $K_{k+2}$, so these children form the
required level $\mathcal T_{k+1}$. The root has descendants at arbitrarily large levels. Since each node has finitely many children, every node with this property has a child
with the same property. Thus there is an infinite branch
$(v_k)_{k\in\mathbb N_0}$ with $v_k\in\mathcal T_k$. Put
\[
F_i^{(k)}=F_i(v_k)\qquad(0\leq i\leq k).
\]
Along this branch, $F_i^{(k)}$ is nonempty and closed in $K_i$ and
\[
\begin{aligned}
&F_i^{(k+1)}\subseteq F_i^{(k)}
   &&(0\leq i\leq k),\\
&F_i^{(k)}\cap\bigcap_{b\in B}\downarrow b\neq\varnothing
   &&(0\leq i<k,\ B\subseteq F_{i+1}^{(k)},\
      1\leq |B|\leq k).
\end{aligned}
\]
Set $E_i=\bigcap_{k\geq i}F_i^{(k)}$. By Lemma~5.9, one can choose
$p_i\in E_i\subseteq K_i$ below every point of $E_{i+1}$. In
particular, $p_i\leq p_{i+1}$ for every $i\in\mathbb N_0$.
\end{proof}

\begin{lemma}
Let $F$ be a nonempty Scott-closed subset of a dcpo $P$. Then $F$
with the induced order is a dcpo and $\sigma(F)=\{F\cap U:U\in\sigma(P)\}$.
If $\Sigma P$ is locally compact, then $\Sigma F$ is locally compact.
If $P$ is meet-continuous, then $F$ is meet-continuous.
\end{lemma}

\begin{proof}
Let $D\subseteq F$ be directed. Since $F$ is Scott closed in $P$,
$\bigvee_P D\in F$. For every $y\in F$,
$y$ is an upper bound if and only if $\bigvee_P D\leq y$. Hence $F$ is a dcpo and $\bigvee_F D=\bigvee_P D$.

Let $C\in\Gamma(F)$. Since $C$ is lower in $F$ and $F$ is lower
in $P$, the set $C$ is lower in $P$. For every directed $D\subseteq C$, $\bigvee_P D=\bigvee_F D\in C$. Thus $C\in\Gamma(P)$. Conversely, if $C\in\Gamma(P)$, then $F\cap C$
is lower in $F$ and every directed $D\subseteq F\cap C$ satisfies
$\bigvee_F D=\bigvee_P D\in F\cap C$. Consequently,
\[
\Gamma(F)=\{F\cap C:C\in\Gamma(P)\}\mbox{ and }
\sigma(F)=\{F\cap U:U\in\sigma(P)\}.
\]

Suppose that $\Sigma P$ is locally compact. Let $x\in F\cap U$ with
$U\in\sigma(P)$. Choose $V\in\sigma(P)$ and a compact saturated
subset $K$ of $\Sigma P$ such that $x\in V\subseteq K\subseteq U$. Since $F$ is closed in $\Sigma P$, the set $F\cap K$ is closed in
$K$ hence compact. Thus we know
\[
x\in F\cap V\subseteq F\cap K\subseteq F\cap U,F\cap V\in\sigma(F),
\]
which proves local compactness of $\Sigma F$.

Finally, suppose that $P$ is meet-continuous. Let $D\subseteq F$ be
directed and let $x\in F$ satisfy $x\leq\bigvee_F D=\bigvee_P D$.
Since $F$ is lower in $P$,
\[
\dn_P x=\dn_F x\mbox{ and }
\dn_P D=\dn_F D.
\]
By meet-continuity of $P$ and the subspace topology on $F$,
\[
x\in F\cap
\operatorname{cl}_{\Sigma P}
   (\dn_P x\cap\dn_P D)
=
\operatorname{cl}_{\Sigma F}
   (\dn_F x\cap\dn_F D).
\]
Hence $F$ is meet-continuous.
\end{proof}

\begin{proposition}
{\rm Let $P$ be a meet-continuous locally compact dcpo and $F$ a nonempty irreducible Scott closed subset of $P$. For every sequence $(O_n)_{n\in\mathbb N_0}$ of nonempty
relatively open subsets of $F$, there exists an increasing sequence
$(p_n)_{n\in\mathbb N_0}$ in $F$ with $p_n\in O_n$ for every $n$.
In particular,
\[
\bigvee_{n\in\mathbb N_0}p_n\in\bigcap_{n\in\mathbb N_0}O_n.
\]}
\end{proposition}

\begin{proof}
By Lemma~5.11, $F$ is a meet-continuous locally compact dcpo and each $O_n$ belongs to $\sigma(F)$. Choose $x_0\in O_0$.
As $\Sigma F$ is locally compact, there exists a compact saturated subset $K_0$
of $\Sigma F$ such that
\[
x_0\in\operatorname{int}_F K_0\subseteq K_0\subseteq O_0.
\]

Suppose that $K_n\subseteq O_n$ has been chosen with
$\operatorname{int}_F K_n\neq\varnothing$. Since $\Sigma F$ is
irreducible and both $O_{n+1}$ and $\operatorname{int}_F K_n$ are
nonempty open subsets of $\Sigma F$, we can choose
\[
x_{n+1}\in O_{n+1}\cap\operatorname{int}_F K_n.
\]
By the local compactness again, there exists a compact saturated subset
$K_{n+1}$ of $\Sigma F$ such that
\[
x_{n+1}\in\operatorname{int}_F K_{n+1}
\subseteq K_{n+1}
\subseteq O_{n+1}\cap\operatorname{int}_F K_n.
\]
Thus each $K_n$ is nonempty and compact in $F$. Since $K_n$
is saturated, we have
\[
K_n\subseteq O_n\mbox{ and }
K_{n+1}\subseteq\operatorname{int}_F K_n
=\operatorname{int}_F(\uparrow_F K_n)=h_F(K_n)
(n\in\mathbb N_0).
\]

By Theorem~5.10, we can choose
$p_n\in K_n$ with $p_n\leq p_{n+1}$ for every $n$. This sequence is
directed. So $p=\bigvee_{n\in\mathbb N_0}p_n$ exists in $F$. For
every $n\in\mathbb N_0$, As $O_n$ is upper in $F$, we obtain
\[
p_n\in K_n\subseteq O_n\mbox{ and }
p_n\leq p\Longrightarrow p\in O_n.
\]
Therefore $p\in\bigcap\limits_{n\in\mathbb N_0}O_n$.
\end{proof}

\begin{corollary}
{\rm Let $P$ be a meet-continuous dcpo. If $\Sigma P$ is locally compact
and first countable, then $\Sigma P$ is sober.}
\end{corollary}

\begin{proof}
Let $F$ be an irreducible Scott closed subset of $P$ and fix
$a,b\in F$. We restrict countable open neighbourhood bases at $a$ and $b$
in $\Sigma P$ to $F$ and denote the resulting bases by
$(U_n)_{n\in\mathbb N_0}$ and $(V_n)_{n\in\mathbb N_0}$, respectively.
By irreducibility of $F$, each $O_n=U_n\cap V_n$ is a nonempty
relatively open subset of $F$.By Proposition~5.12, it gives an increasing
sequence $(p_n)_{n\in\mathbb N_0}$ in $F$ such that
$p_n\in O_n$ and
$$c=\bigvee_{n\in\mathbb N_0}p_n
\in\bigcap_{n\in\mathbb N_0}O_n
=
\Bigl(\bigcap_{n\in\mathbb N_0}U_n\Bigr)
\cap
\Bigl(\bigcap_{n\in\mathbb N_0}V_n\Bigr).
$$

We have $\bigcap\limits_{n\in\mathbb N_0}U_n=\up_F a$. Indeed, each
$U_n$ is upper in $F$ and contains $a$. So $\up_F a\subseteq U_n$.
Conversely, if $z\in F\setminus\up_F a$, then $a\not\leq z$ and
\[
a\in F\setminus\dn_F z
  =F\cap(P\setminus\dn_P z).
\]
Since $\dn_P z$ is Scott closed, there exists $n\in\mathbb N_0$
such that $a\in U_n\subseteq F\setminus\dn_F z$. Hence
$z\notin U_n$. The same argument at $b$ gives
$\bigcap_{n\in\mathbb N_0}V_n=\up_F b$. Therefore $c\in\up_F a\cap\up_F b$. Thus $F$ is directed. Put $p=\bigvee_P F$. Since $F$ is Scott closed,
$p\in F$. So we have $F=\dn_P p=\operatorname{cl}_{\Sigma P}\{p\}$. Hence $\Sigma P$ is sober.
\end{proof}

\begin{lemma}
{\rm Let $P$ be a meet-continuous locally compact dcpo and $F$ an irreducible Scott closed subset
of $P$. For every sequence $(H_n)_{n\in\mathbb N_0}$ of nonempty
relatively open subsets of $F$, there exists
$b\in\bigcap_{n\in\mathbb N_0}H_n$ having a countable open
neighbourhood base in $F$. In particular, $\up_F b$ is a
$G_\delta$ subset of $F$.}
\end{lemma}

\begin{proof}
By Lemma~5.11, $F$ is meet-continuous and $\Sigma F$ is locally
compact with the subspace topology. Put $W_{-1}=F$. Given a nonempty
$W_{n-1}\in\sigma(F)$, the irreducibility gives
$W_{n-1}\cap H_n\neq\varnothing$. By local compactness, choose a
nonempty $W_n\in\sigma(F)$ and a compact saturated subset $Q_n$ of
$\Sigma F$ such that
\[
W_n\subseteq Q_n\subseteq W_{n-1}\cap H_n.
\]
Thus for every $n\geq0$, $Q_{n+1}\subseteq W_n\subseteq Q_n\subseteq H_n$.
Consider the family
\[
\mathcal C=
\bigl\{B\in\Gamma(F):
       B\cap Q_n\neq\varnothing\text{ for every }n\geq0\bigr\}.
\]
We have $F\in\mathcal C$. Let $\mathcal A$ be a nonempty chain in
$\mathcal C$ under inclusion. For each $n\geq0$, the sets
$A\cap Q_n$ ($A\in\mathcal A$) are closed in $Q_n$. For every
nonempty finite $\mathcal A_0\subseteq\mathcal A$,
\[
\bigcap_{A\in\mathcal A_0}(A\cap Q_n)
=(\min\nolimits_{\subseteq}\mathcal A_0)\cap Q_n
\neq\varnothing.
\]
By compactness of $Q_n$, $Q_n\cap\bigcap_{A\in\mathcal A}A\neq\varnothing
(n\geq0)$.
Hence $\bigcap\mathcal A\in\mathcal C$. By Zorn's lemma,
$\mathcal C$ has a minimal member $B$.

Let $O\in\sigma(F)$ with $O\cap B\neq\varnothing$. Since
$B\setminus O$ is closed in $F$ and properly contained in $B$,
the minimality gives
\[
\exists j\geq0,
(B\setminus O)\cap Q_j=\varnothing,
\mbox{ and}
\varnothing\neq B\cap Q_j\subseteq O.
\]
If $O_1,O_2\in\sigma(F)$ both meet $B$, choose corresponding
indices $j_1,j_2$ and put $j=\max\{j_1,j_2\}$. Then
\[
\varnothing\neq B\cap Q_j
\subseteq (B\cap Q_{j_1})\cap(B\cap Q_{j_2})
\subseteq O_1\cap O_2.
\]
Thus $B$ is irreducible. Since $B$ is closed in $F$ and $F$ is
Scott closed in $P$, $B$ is also Scott closed in $P$.

For each $n\geq0$, the set $B\cap W_n$ is relatively open in $B$ and $\varnothing\neq B\cap Q_{n+1}\subseteq B\cap W_n$.

By Proposition~5.12, we have $b\in B\cap\bigcap\limits_{n\geq0}W_n$. In particular, $b\in Q_n\cap H_n$ for every $n$. Since $B$ is lower
in $F$ and $b\in B$, we have $\downarrow_F b\subseteq B$. Moreover,
$\dn_F b$ is closed in $F$ and $b\in(\dn_F b)\cap Q_n(n\geq0)$.
Hence $\dn_F b\in\mathcal C$ and the minimality yields
$B=\dn_F b$.

For each $n\geq0$, put
\[
V_n=\up_F(W_n\cap\dn_F b).
\]
Since $b\in W_n\in\sigma(F)$, by Lemma~5.1,
$b\in V_n\in\sigma(F)$. Let $O\in\sigma(F)$ with $b\in O$.
As $b\in B\cap O$, choose $j\geq0$ with $B\cap Q_j\subseteq O$.
Then $W_j\cap\downarrow_F b=W_j\cap B\subseteq Q_j\cap B\subseteq O$.
Since $O$ is upper in $F$, it follows that $b\in V_j\subseteq O$.
Thus $(V_n)_{n\geq0}$ is a countable open neighbourhood base at $b$. Finally, $\uparrow_F b\subseteq V_n$ for every $n$. If
$z\in F\setminus\up_F b$, then $F\setminus\dn_F z$
is an open neighbourhood of $b$, so some $V_j$ is contained in
$F\setminus\dn_F z$. Hence $z\notin V_j$. This proves $\up_F b=\bigcap\limits_{n\geq0}V_n$.
\end{proof}

\begin{theorem}
Let $P$ be a meet-continuous locally compact dcpo. Then $\Sigma P$ is sober.
\end{theorem}

\begin{proof}
Let $F$ be a nonempty irreducible Scott closed subset of $P$ and put
\[
D=\{x\in F:x\text{ has a countable open neighbourhood base in }F\}.
\]
By Lemma~5.11, $\sigma(F)$ is the subspace topology on $F$.
For each nonempty $U\in\sigma(F)$, using Lemma~5.14 to the constant
sequence $H_n=U$, $n\in\mathbb N_0$, we have
$D\cap U=D\cap\bigcap_{n\geq0}H_n\neq\varnothing$.
Hence $D\neq\varnothing$ and $\operatorname{cl}_{\Sigma F}D=F$.
For $x\in D$ and any countable open neighbourhood base
$(U_n)_{n\geq0}$ at $x$ in $F$, we have
\[
\bigcap_{n\geq0}U_n=\up_F x.
\]
Indeed, each $U_n$ is upper in $F$ and contains $x$. Conversely,
if $z\in F\setminus\up_F x$, then
\[
x\in F\setminus\dn_F z
 =F\cap(P\setminus\dn_P z)\in\sigma(F).
\]
Thus $U_n\subseteq F\setminus\dn_F z$ for some $n$ and
$z\notin U_n$.

Now fix $a,b\in D$ and choose countable open neighbourhood bases
$(U_n)_{n\geq0}$ at $a$ and $(V_n)_{n\geq0}$ at $b$ in $F$.
Put $H_{2n}=U_n$ and $H_{2n+1}=V_n$. By Lemma~5.14, there exists
\[
c\in D\cap\bigcap_{n\geq0}H_n
 =D\cap\Bigl(\bigcap_{n\geq0}U_n\Bigr)
    \cap\Bigl(\bigcap_{n\geq0}V_n\Bigr)
 =D\cap\up_F a\cap\up_F b.
\]
Hence $a\leq c$ and $b\leq c$. So $D$ is directed. Put $p=\bigvee_P D$. Since $D\subseteq F$ and $F$ is Scott closed,
$p\in F$. The set $F\cap\dn_P p$ is closed in $\Sigma F$ and
contains $D$. Since $F$ is lower, it follows that
\[
F=\operatorname{cl}_{\Sigma F}D
 \subseteq F\cap\dn_P p
 \subseteq\dn_P p
 \subseteq F.
\]
Thus $F=\dn_P p=\operatorname{cl}_{\Sigma P}\{p\}$. Therefore $\Sigma P$ is sober.
\end{proof}

\vspace{5pt}

{\bf AI Declaration} For more readable, the authors just use AI tools to assist with English language editing. The author independently verified the mathematical statements, proofs and references. The author takes full responsibility for the manuscript.

\end{document}